\documentstyle[11pt]{amsart}

\newtheorem{theorem}{Theorem}
\newtheorem{lemma}{Lemma}
\newtheorem{corollary}{Corollary}
\theoremstyle{remark}
\newtheorem{remark}{Remark} 

\newcommand{\thmfont}{}

\newcommand{\complex}{{\mathbb C}}

\newcommand{\ox}{o_{\text{\rm excl}}}
\newcommand{\Ox}{O_{\text{\rm excl}}}
\def\spa#1{\complex^#1}
\renewcommand{\cal}{\mathcal}
\newcommand{\M}[1]{M^{(#1)}}
\newcommand{\N}{{\cal N}}
\newcommand{\m}[1]{m^{(#1)}}
\newcommand{\A}[1]{a^{(#1)}}
\newcommand{\p}[2]{\frac{\partial #1}{\partial #2}}
\renewcommand{\a}{\alpha}
\newcommand{\Ker}{\operatorname{Ker}}

\begin{document}

\title{Ordinary differential equations with only entire solutions}

\author{Erik Anders\'en}


\address{Department of Mathematics\\ Van Vleck Hall\\ 480 Lincoln
Drive\\ University of Wisconsin-Madison\\ 
WI 53706}

\email{andersen@@math.wisc.edu}

\thanks{Supported by a grant from STINT}

\begin{abstract}
We prove necessary and sufficient conditions for a system $\dot z_i=z_ip_i(z)$
($p_i$ a polynomial) to have only entire analytic functions as solutions.
\end{abstract}

\maketitle

\section{preliminaries}
We will need the following elementary facts from Nevanlinna theory. For any
entire function of one complex variable $f(z)$ we set
\begin{gather*}
m(r,f)=\frac1{2\pi r}\int_{|z|=r}\log^+|f(z)|\,|dz|,\\
N(r,f)=\sum_{f(a)=\infty}\log^+\left|\frac ra\right|\quad\text{and}\\
T(r,f)=m(r,f)+N(r,f).
\end{gather*}
In the definition of $N$, the sum is taken over all poles $a$ of $f$, with
regard to multiplicity. The function $T(r,f)$ is called the Nevanlinna 
characteristic of $f$. The following properties
are easily verified.
\begin{equation} \label{Nevanlinna-rules}
\begin{gathered} 
T(r,f_1+f_2)\le T(r,f_1)+T(r,f_2)+O(1)\\
T(r,f_1f_2)\le T(r,f_1)+T(r,f_2),\\
T(r,f^d)=dT(r,f).
\end{gathered}
\end{equation}
Here and in the sequel, the estimates $O(g)$ and $o(g)$ are as $r\to\infty$.
The first fundamental theorem of Nevanlinna theory says that
\[
T(r,1/f)=T(r,f)+O(1).
\]
We also have the Lemma of the logarithmic derivative (LLD),
\[
m(r,f'/f)=\ox(T(r,f)).
\]
Here $\ox$ means that the estimate holds outside a set
of finite measure. We also use the corresponding notation $\Ox$. From LLD
and the preceding inequalities it follows that
\begin{gather}
m(r,f^{(k}))\le(1+\ox(1)))T(r,f)\\
m(r,f^{(k)}/f)=\ox(T(r,f))
\end{gather}
for all positive integers $k$.

\section{Borel's Theorem} We need a version of a theorem named after Borel,
which is classical, but which is  proved here since I have not found a
good reference to it. We call an entire holomorphic function a unit if
it has no zeros. 

\begin{theorem}\thmfont Let $f$ be an entire function and $u_1$, $u_2$, $\dots$,
$u_n$ be units satisfying  
\begin{equation}\label{u-sum}
\sum u_i=f.
\end{equation}
Then one of the following cases holds.
\begin{enumerate}
\item $T(r,u_i)=\Ox(T(r,f)+1)$ for all $i$.
\item Some subsum $\sum_{i\in I}u_{i}=0$.
\end{enumerate}
\end{theorem}

\begin{pf} The proof is by induction. If $n=1$ then Case~1 holds automatically
 so there is nothing to prove. We assume that the theorem holds for sums with 
less  than $n$ terms. If we differentiate (\ref{u-sum})  we get
\[
\sum \frac{u_i^{(j)}}{u_i}u_i=f^{(j)}.\quad\text{for $j=0,\dots,n-1$}
\]
This is a linear system for $u_i$. Two cases are possible.
\begin{enumerate}
\item $\det(u_i^{(j)}/u_i)\not\equiv0$. Then we can solve for $u_i$ 
and get
\begin{multline*}
T(r,u_i)=O(T(r,f))+\sum_j O(T(r,u_i^{(j)}/u_i))+O(1)\\
=O(T(r,f))+\sum_i\ox(T(r,u_i))+O(1)
\end{multline*}
by LLD. We therefore have Case~1.
\item $\det(u_i^{(j)}/u_i)\equiv0$. This means that the Wronskian
$W(\{u_i\})=\det(u_i^j)\equiv0$.
The theory of ordinary differential equations now says there are constants
$c_i$ such that
\begin{equation}\label{u dependence}
\sum c_iu_i=0.
\end{equation}
We choose the shortest possible such sum, that is, the sum with the smallest
number of non-zero $c_i$.
Possibly after a reordering and a scaling
we may assume that $c_1=-1$, so that $u_1=\sum_{i>1} c_iu_i$. Division by
$u_1$ gives
\begin{equation}\label{v-sum}
\sum_{i>1}c_iu_i/u_1=1.
\end{equation}
We set $v_i=c_iu_i/u_1$ and apply the theorem to the expression (\ref{v-sum}),
which has less than $n$ terms. There are two cases.
\begin{enumerate}
\item $T(v_i,r)=\Ox(T(1,r)+1)=\Ox(1)$ for all $i$. Then all $v_i$ are constant
and $u_i=a_iu_1$ for some $a_i$. The sum (\ref{u-sum}) now becomes
$(\sum a_i+1)u_1=f$. If $\sum a_i+1=0$, then $f=0$ and we have Case~2;
otherwise, $u_1$ is proportional to $f$ and then all $u_i$ are also
proportional to $f$. We then have Case~1.
\item $\sum_{i\in I}v_i=0$ for some set $I$. Then $\sum_{i\in I}c_iu_i=0$ and
this sum is shorter than (\ref{u dependence}). This is a contradiction. \qed

\end{enumerate}

\end{enumerate}
\renewcommand{\qed}{}
\end{pf}

\section{Differential equations}

\begin{theorem}\label{ode-theorem}\thmfont Let $p_i$ be Laurent-polynomials and
 \begin{equation}\label{ode}
\dot z_i=z_ip_i(z)\qquad z_i(0)=c_i\qquad 1\le i\le n
 \end{equation}
be an ordinary differential equation such that for all $c$ in some open
set all components of the solutions are units. Then
 \begin{equation}\label{result}
p(Z)=\sum_{k=1}^{n-1} u_k\theta_k(M^{(k)})+ u_0.
 \end{equation}
Here $u_k\in\spa{n}$, $\M k\in\complex[Z,Z^{-1}]^k$ is defined by $\M n=Z$ and
\begin{equation}\label{M-equation}
\M{k-1}_i=\prod{\M{k}_j}^{\A k_{ij}}\text{for $k=n,\dots,2$},
\end{equation}
where $\A
k$ are $(k-1,k)$-matrices satisfying $\A k\dots\A nu_k=0$
and $\theta_k$ are Laurent-polynomials. Conversely, any system of the
form~(\ref{result}) has only entire solutions, and all components are either
units or identically zero.
\end{theorem}

The proof needs a lemma. We write $p(Z)=(p_1(Z),\dots,p_n(Z))$,
$p(Z)=\sum p_\a Z^\a$ and $p_\a=(p_{\a,1},\dots,p_{\a,n})$. We let $Dg$ denote
$(\dot g_1/g_1,\dots,\dot g_n/g_n)$. Let $S$ be the set of meromorphic
functions satisfying $T(r,f)=\ox(T(r,z_1)+\dots+T(r,z_n)+1)$.

\begin{lemma}\label{dimension}\thmfont Let ${\cal M}$ be the
multiplicative module generated by of
${\cal N}=\{Z^\a:p_{\a,i}\ne0\text{ for some }i\}$. Then $\dim{\cal M}<n$.
\end{lemma}

\begin{pf}  Choose $c$ so that $\sum_{\a\in A} p_{\a,i} c^\a\ne0$ for
all non-zero subpolynomials and all indices $i$, and so that all
components of the  solution are units. Let $z(t)$ be the
solution with $z(0)=c$. Apply Borel's theorem to $\sum p_{\a,i}z^\a=
\dot z_i/z_i$. Because of our choice of initial condition, Case~2 does
not hold. Therefore, $\N\subset S$. It follows from
the rules (\ref{Nevanlinna-rules}) that ${\cal M}\subset S$. Assume that the
dimension of $\cal M$ is~$n$. Then for each $i$ there is $d_i$ such that
$z^{d_i}\in{\cal M}\subset S$. Again by
(\ref{Nevanlinna-rules}), $z_i\in S$ for each $i$. This is impossible
and therefore the dimension is less than $n$. 
\end{pf}

\begin{pf*}{Proof of Theorem \ref{ode-theorem}} The proof
is by induction. The theorem
holds if $n=1$, so we assume that  it holds for systems with smaller
number of variables. We first prove the necessity of (\ref{result}).
Assume that (\ref{ode}) has only entire solutions. By
Lemma~\ref{dimension} we can 
write $p(Z)=f(M)$, where $M=(M_1,\dots,M_{l})$ and $M_i=\prod Z_j^{a_{ij}}$,
and $l<n$. We may assume that $a=(a_{ij})$ has full rank since we can
otherwise drop some variables $M_i$. If $l<n-1$ we can complete $a$ to an
$(n-1,n)$-matrix of full rank, so we assume that $l=n-1$. Set
$m_i=\prod z_j^{a_{ij}}$. Then $Dz=p(z)=f(m)$ and 
logarithmic differentiation gives $Dm=af(m)$. By induction,
 \[
af(M)=\sum_{k=1}^{n-2}v_k\theta_k(\M k)+v_0,
 \]
where $\M k$ are given by (\ref{M-equation}) for $k=n-1,\dots,2$ and
$\A k\dots\A{n-1}v_k=0$.
Since $a$ is surjective $v_k=au_k$ for some $u_k$ so
 \[
f(M)-\sum_{k=1}^{n-2}v_k\theta_k(\M k)-v_0\in\Ker a.
 \]
Since this kernel is  one-dimensional we have
 \[
f(M)=u_{n-1}\theta_{n-1}(M)+\sum_{k=1}^{n-2}u_k\theta_k(\M k)+u_0
 \]
for some $u_{n-1}\in\Ker a\subset\spa{n}$ and some Laurent-polynomial
$\theta_{n-1}$.
If we set $\M{n-1}=M$ and $\A{n-1}=a$ the first half of the theorem is
proved.

If the  system is  of form (\ref{result}) then we get as above that
 \begin{align*}
Dm&=\A n(\sum_{k=1}^{n-1}u_k\theta_k(\m k)+u_0)\\
  &=\sum_{k=1}^{n-2}\A nu_k\theta_k(\m k)+\A nu_0
 \end{align*}
since $\A nu_n=0$ by assumption. By induction, this system has only
 entire  solutions $m(t)$. An integration of $\dot z_i/z_i=f_i(m(t))$
 now proves the rest of the theorem.
\end{pf*}

\begin{remark} Although the form (\ref{result}),
(\ref{M-equation}) is complicated, the 
proof gives an extremely simple recursive algorithm to determine if
all solutions 
to a system are entire. First compute the multiplicative span of all
$Z^\a$ such that
$p_{\a,i}\ne0$ for some $i$. If this is the whole space, the system does
have non-entire solutions. Otherwise, let $M_i$ be a basis for this space,
and express $p(Z)=f(M)$ in terms of this basis. Differentiate logarithmically
to get a new equation $Dm=af(m)$ (with notations as in the theorem). The
original system has only entire solutions if and only if this new system has.
\end{remark}
\begin{remark} If the polynomials $p_i$ are genuine polynomials, it is
immediately clear that all components of any solution are either zero-free
or identically zero.
\end{remark}
\begin{corollary}\thmfont All complete polynomial vector fields on
$(\complex^{*})^n$ are of 
the form (\ref{ode}) with $p$ satisfying (\ref{result}).
\end{corollary}
\begin{corollary}\thmfont All complete polynomial vector fields on
$(\complex^{*})^n$ preserve 
the volume form $\bigwedge dz_i/z_i$.
\end{corollary}
\begin{pf} The proof is again by induction. We compute
 \[
\frac{d}{dt}\frac{dz_i}{z_i}=\frac{d\dot z_i\,z_i-\dot z_i\,dz_i}{z_i^2}
=\frac{d(z_ip_i(z))z_i-z_ip_i(z)dz_i}{z_i^2}=d(p_i(z)).
\]
This shows in particular that the result holds for $n=1$. Also, we compute
\[
 \frac{d}{dt}\bigwedge_i \frac{dz_i}{z_i}=
\sum_i\frac{dz_1}{z_1}\wedge\dots\wedge\frac{d}{dt}\frac{dz_i}{z_i}\wedge
\dots\wedge\frac{dz_n}{z_n}=\left(\sum_i z_i\p{p_i}{z_i}\right)
\left(\bigwedge_i\frac{dz_i}{z_i}\right).
\]
We use the notation of the proof of the theorem and in particular the result
that $p(z)=f(m)$ and $Dm=af(m)$. We have to prove that
$\sum z_i{\partial p_i}/{\partial z_i}(z)=0$,
so we compute
\begin{multline}
\sum_i z_i\p{f_i}{z_i}=\sum_i z_i\sum_j\p{f_i}{m_j}\p{m_j}{z_i}=
\sum_i z_i\sum_j\p{f_i}{m_j}m_j\frac{a_{ji}}{z_i}\\
=\sum_{i,j}\p{f_i}{m_j}m_ja_{ji}=\sum_j\p{(af)_j}{m_j}m_j.
\end{multline}
The last expression is $0$ by induction and the corollary is proved.
\end{pf}
\end{document}